\documentclass[letterpaper,10pt]{article}

\usepackage[utf8]{inputenc}
\usepackage[english,brazil]{babel}
\usepackage[T1]{fontenc}

\usepackage{comment}

\usepackage{amsmath, amsfonts, amssymb, amsthm}

\theoremstyle{plain}

\newtheorem{definition}{Definição}
\theoremstyle{definition}
\newtheorem{proposition}{Proposição}

\theoremstyle{remark}

\usepackage{graphicx}

\usepackage[hidelinks]{hyperref}

\usepackage[left=3.4cm,right=3.4cm,bottom=3.5cm,top=3.4cm]{geometry}

\title{Fundamentos para a Construção no GeoGebra de Tesselações Aperiódicas usando um Único Polígono}
\author{Astor Santos Neto, Sandra Maria Barbosa, Alcebiades Dal Col\footnote{dilemastor@gmail.com, sms\underline{ }barbosa@hotmail.com, alcebiades.col@ufes.br} \\ \\ Universidade Federal do Espírito Santo}

\begin{document}

\maketitle

\begin{abstract}
Neste trabalho, nós iremos explorar alguns polígonos que individualmente são capazes de preencher o plano de forma aperiódica. Estes polígonos foram descobertos recentemente por alguns pesquisadores e configuram um grande achado para a Matemática. Nós apresentaremos duas formas de construção para estes polígonos, a primeira delas é fundamentada na justaposição de copias de uma pipa em particular e a segunda consiste em construir, com régua e compasso, a linha poligonal formada pelos lados de cada um destes polígonos. As construções são feitas no GeoGebra com riqueza de detalhes para permitir a reprodução.

\medskip

\noindent \textbf{Palavras-chave:} Tesselações aperiódicas; Chapéu; Tartaruga; GeoGebra
\end{abstract}

\begin{otherlanguage}{english}
\begin{abstract}
In this work, we will explore some polygons that individually are capable of filling the plane in an aperiodic way. These polygons were recently discovered by some researchers and constitute a great discovery for Mathematics. We will present two ways of constructing these polygons, the first of which is based on the juxtaposition of copies of a particular kite and the second consists of constructing, with a ruler and compass, the polygonal line formed by the sides of each of these polygons. The constructions are made in GeoGebra with great detail to allow reproduction.

\medskip

\end{abstract}
\end{otherlanguage}
 
\section{Introdução}

Os ladrilhos estão naturalmente no nosso cotidiano em caminhos e paredes. Por exemplo, podemos ter uma parede preenchida apenas com hexágonos $(6. 6. 6)$, um chão preenchido com triângulos e dodecágonos $(3. 12. 12)$, ou um teto formado com quadrados e octógonos $(4. 8. 8)$. Os ladrilhamentos também podem ser chamados de mosaicos ou tesselações.

No mundo matemático, dizemos que um conjunto de polígonos gera uma tesselação do plano se este conjunto for capaz de preencher completamente uma superfície bidimensional sem a sobreposição de polígonos e sem a existência de espaços vazios. As tesselações podem ser classificadas como periódicas ou aperiódicas~\cite{neto2023penroseGeoGebra}. Dentre as tesselações aperiódicas existentes, nós destacamos as tesselações descobertas recentemente por Smith e seus colaboradores~\cite{smith2023aperiodic}, cujo trabalho deu continuidade as investigações a respeito de conseguir uma tesselação aperiódica com o menor número de polígonos possível~\cite{neto2023penroseGeoGebra}. Smith e seus colaboradores descobriram polígonos que individualmente são capazes de preencher completamente o plano de forma aperiódica.

Já no mundo educacional, o estudo de ladrilhamentos é recomendado pela Base Nacional Comum Curricular (BNCC)~\cite{bncc2018bme}, como vemos nas competências (EF07MA27) e (EM13MAT505). É sugerido que este estudo pode ser feito com o uso de aplicativos de geometria dinâmica que auxiliam o aluno em seu aprendizado, por exemplo: o \emph{software GeoGebra}. Logo, com o objetivo de auxiliar alunos e professores, alguns conceitos e construções de tesselações aperiódicas serão expostos nas próximas seções.

\begin{itemize}
    \item (EF07MA27) Calcular medidas de ângulos internos de polígonos regulares, sem o uso de fórmulas, e estabelecer relações entre ângulos internos e externos de polígonos, preferencialmente vinculadas à construção de mosaicos e de ladrilhamentos. \cite[p.~309]{bncc2018bme}
    \item (EM13MAT505) Resolver problemas sobre ladrilhamento do plano, com ou sem apoio de aplicativos de geometria dinâmica, para conjecturar a respeito dos tipos ou composição de polígonos que podem ser utilizados em ladrilhamento, generalizando padrões observados. \cite[p.~541]{bncc2018bme}
\end{itemize}

\section{Tesselações Semirregulares e Tesselações Duais de Laves}

As tesselações podem ser classificadas como regulares, semirregulares, demirregulares ou irregulares. As tesselações regulares são feitas com repetições de um mesmo polígono regular, de modo que existem três possíveis tesselações: com triângulos, quadrados e hexágonos. Isso acontece porque para revestir o plano utilizando apenas um tipo de polígono regular sem deixar lacunas e sem haver sobreposição, os ângulos internos desse polígono devem ser divisores exatos de $360^{\circ}$~\cite{CAVALCANTI}. Neste trabalho, iremos nos concentrar nas tesselações semirregulares, estas são feitas de dois ou mais polígonos regulares diferentes. Um exemplo particularmente importante para o nosso estudo é a tesselação $(3.4.6.4)$\footnote{Para mais detalhes sobre esta nomenclatura, veja~\cite{neto2023penroseGeoGebra}.}.

\begin{figure}[ht]
    \centering
    \includegraphics[width=0.7\textwidth]{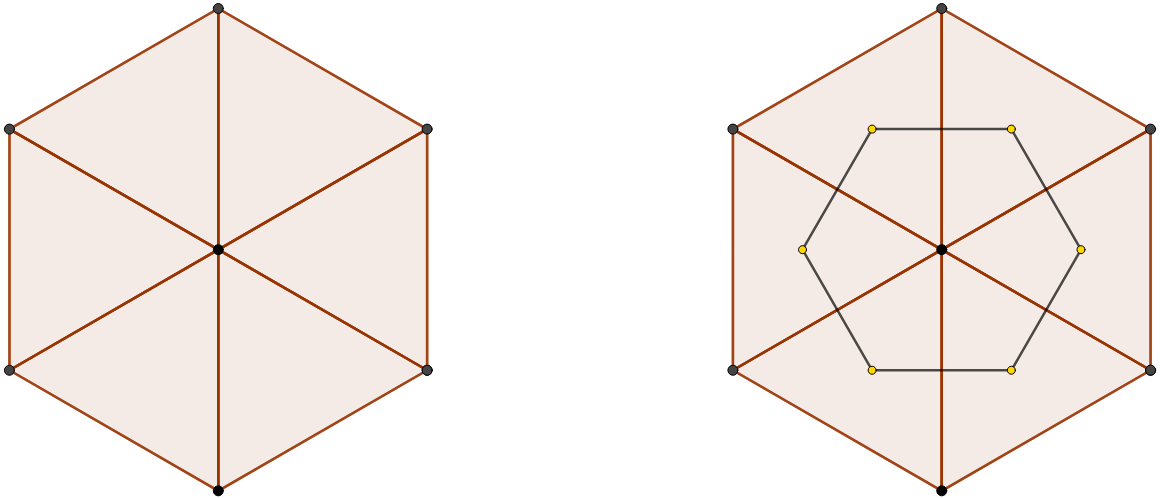}
    \caption{Tesselação dual de Laves para uma tesselação regular com triângulos}
    \label{fig:tesselacao_triangular_Laves}
\end{figure}

Uma \emph{tesselação dual de Laves} é obtida a partir de uma outra tesselação inicial. Mais precisamente, dada uma tesselação, uma tesselação (dual) de Laves têm os seus vértices dados pelos centros dos polígonos da tesselação original e as suas arestas conectam centros cujos polígonos correspondentes tem uma aresta em comum~\cite{grunbaum1987tilings}. Aplicações modernas das tesselações de Laves incluem sistemas de materiais topologicamente interligados~\cite{williams2021mechanics}. Para exemplificar, vamos considerar uma tesselação com triângulos equiláteros $(3.3.3.3.3.3)$. A Figura~\ref{fig:tesselacao_triangular_Laves} mostra à esquerda uma parte de uma tesselação regular com triângulos e à direita os centros dos triângulos equiláteros são exibidos na cor amarela, os quais são conectados caso os triângulos tenham uma aresta em comum. Ao expandirmos esta construção para toda a tesselação com triângulos equiláteros $(3.3.3.3.3.3)$ obtemos uma tesselação dual de Laves que será uma tesselação regular com hexágonos $(6.6.6)$.

\begin{figure}[ht]
    \centering
    \includegraphics[width=0.7\textwidth]{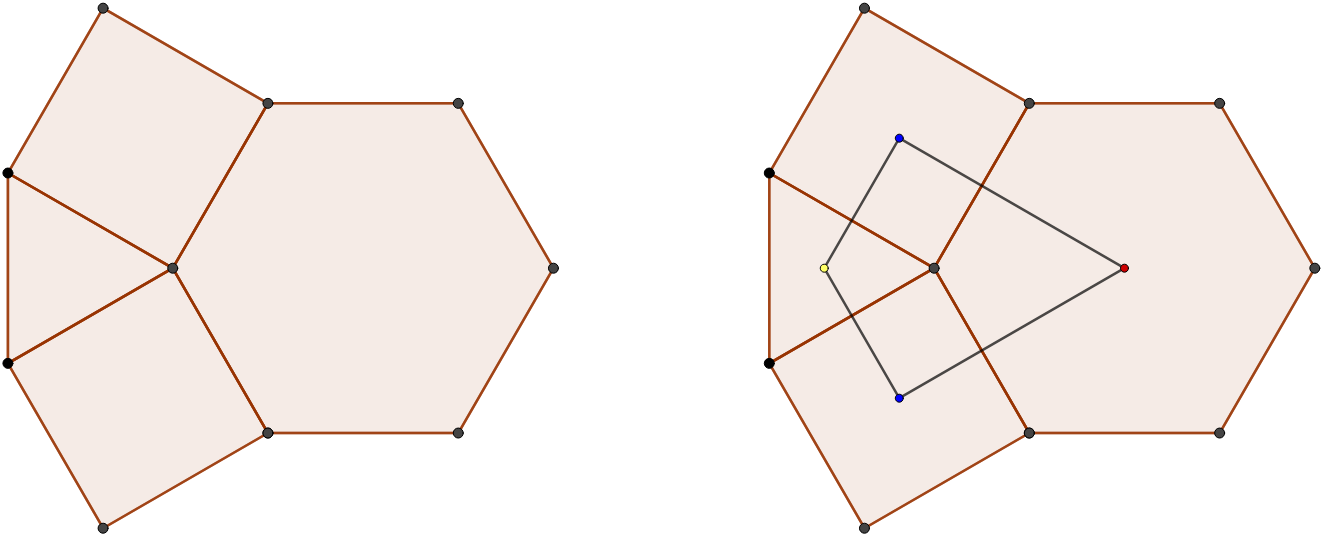}
    \caption{Tesselação dual de Laves para uma tesselação semirregular $(3.4.6.4)$}
    \label{fig:tesselacao_semirregular_Laves}
\end{figure}

Uma segunda tesselação dual de Laves que será importante para os nossos objetivos é a tesselação dual da tesselação semirregular $(3.4.6.4)$. A Figura~\ref{fig:tesselacao_semirregular_Laves} mostra à esquerda uma parte da tesselação semirregular $(3.4.6.4)$ e à direita é exibido a tesselação dual de Laves correspondente. Com esta construção obtemos um objeto importante, um polígono chamado de \emph{pipa}. Mostraremos a seguir que os polígonos de Smith e seus colaboradores podem ser obtidos a partir da justaposição destas pipas seguindo uma determinada ordem de construção. Deste modo, a tesselação dual de Laves da tesselação semirregular $(3.4.6.4)$ é um ambiente propício para a construção destes polígonos.

\section{Construção dos Polígonos de Smith no GeoGebra}

Nesta seção, nós mostraremos inicialmente como construir uma pipa que será fundamental para a construção dos polígonos de Smith e seus colaboradores.

\begin{definition}
\label{def:pipa}
Uma \emph{pipa} é um quadrilátero com dois pares de lados adjacentes congruentes.
\end{definition}

\begin{figure}[ht]
    \centering
    \includegraphics[width=0.4\textwidth]{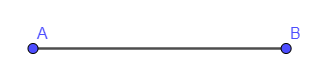}
    \caption{Primeira etapa na construção da pipa}
    \label{fig:pipa_etapa1}
\end{figure}

Iniciamos o software GeoGebra selecionando a \emph{Disposição de Geometria} e configuramos o programa para mostrar a \emph{Janela de Álgebra}. Em seguida, escolhemos dois pontos quaisquer do plano, que serão denotados por $A$ e $B$, e conectamos estes pontos com um segmento. A título de exemplo, usaremos os pontos $A=(0,0)$ e $B=(1,0)$ (veja a Figura~\ref{fig:pipa_etapa1}).

\begin{figure}[ht]
\begin{minipage}[b]{0.48\textwidth}
    \centering
    \includegraphics[width=0.8\textwidth]{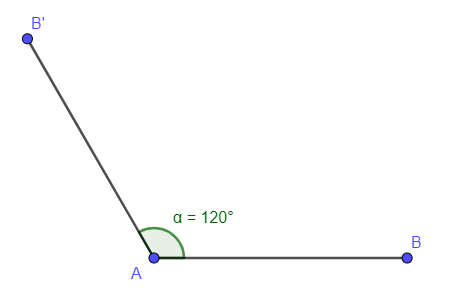}
    \caption{Segunda etapa na construção da pipa}
    \label{fig:pipa_etapa2}
\end{minipage} \quad
\begin{minipage}[b]{0.48\textwidth}
    \centering
    \includegraphics[width=\textwidth]{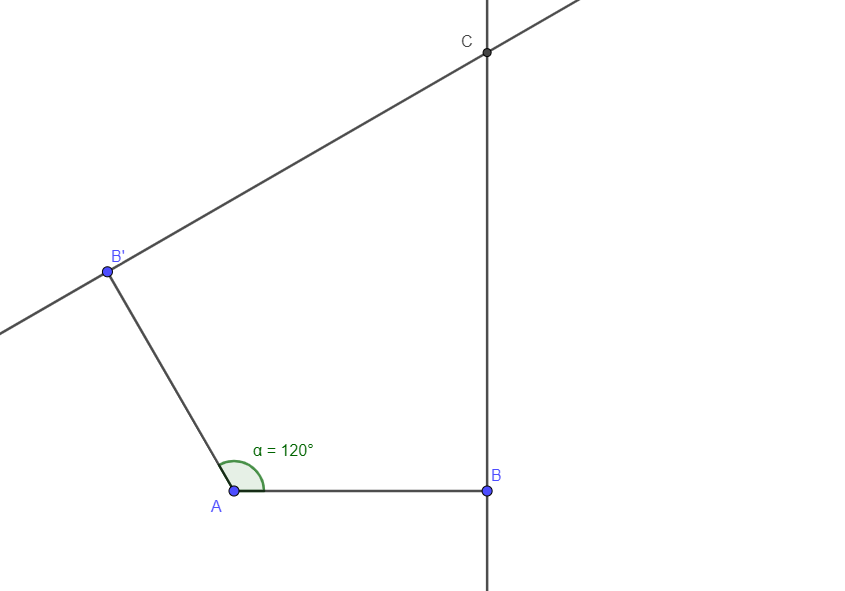}
    \caption{Terceira etapa na construção da pipa}
    \label{fig:pipa_etapa3}
\end{minipage}
\end{figure}

Continuamos a construção criando um ângulo de $120^{\circ}$ no sentido anti-horário a partir dos pontos $B$ e $A$ (nesta ordem) e conectamos o novo ponto $B'$ obtido deste processo com o ponto $A$ usando um segmento de reta (veja a Figura~\ref{fig:pipa_etapa2}).

Neste momento, construímos uma reta tangente ao segmento $AB$ passando pelo ponto $B$ e uma reta tangente ao segmento $AB'$ passando pelo ponto $B'$ e determinamos o ponto $C$ de interseção destas retas tangentes (veja a Figura~\ref{fig:pipa_etapa3}).

\begin{figure}[ht]
    \centering
    \includegraphics[width=0.5\textwidth]{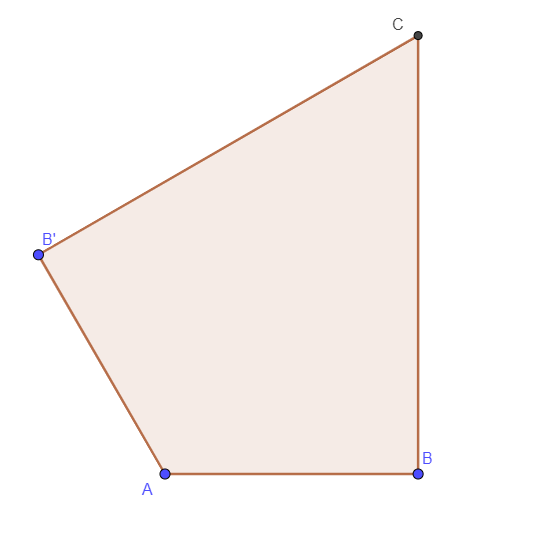}
    \caption{Quarta etapa na construção da pipa}
    \label{fig:pipa_etapa4}
\end{figure}

Finalmente, nós omitimos todos os elementos desta construção com exceção dos pontos (clique Ctrl+A, em seguida, com o botão Ctrl pressionado, clique nos quatro pontos, e omita os objetos ainda selecionados). Em seguida, construímos o polígono $ABCB'$ (Figura~\ref{fig:pipa_etapa4}). 

\begin{proposition}
O polígono $ABCB'$ é uma pipa (veja a Figura~\ref{fig:pipa_etapa4}).
\end{proposition}
\begin{proof}
Ao construirmos a diagonal $AC$, obtemos dois triângulos retângulos com uma hipotenusa em comum e catetos $AB$ e $AB'$ congruentes. Pelo caso especial de congruência para triângulos retângulos, temos que os triângulos $ABC$ e $AB'C$ são congruentes. Consequentemente, os segmentos $BC$ e $B'C$ são congruentes como queríamos verificar para confirmar que $ABCB'$ é uma pipa. Esta pipa será particularmente importante no que segue por este motivo a denominarmos por \emph{pipa de Laves}. Esta nomenclatura foi escolhida visto que este polígono também pode ser obtido a partir da tesselação dual de Laves da tesselação semiregular $(3.4.6.4)$.
\end{proof}

\begin{definition}
Um \emph{$Tile(a,b)$} é um polígono formado pela justaposição de pipas de Laves com um par de lados adjacentes congruentes medindo $a$ e o outro par de lados medindo $b$.
\end{definition}

Passaremos a construção do polígono $Tile(1,\sqrt{3})$ no software GeoGebra. Este polígono ficou mundialmente famoso com o nome de \emph{Hat} e o mesmo pode ser obtido com $8$ pipas de Laves. Em seguida, faremos a construção do polígono $Tile(\sqrt{3},1)$. Este polígono ficou amplamente conhecido como \emph{Turtle} e ele pode ser construído com $10$ pipas de Laves. Finalmente, faremos a construção dos polígonos $Tile(0,1)$, $Tile(1,1)$ e $Tile(1,0)$.
Todas estas construções partem da construção inicial da pipa de Laves apresentada anteriormente.

\begin{figure}[ht]
    \centering
    \includegraphics[width=0.65\textwidth]{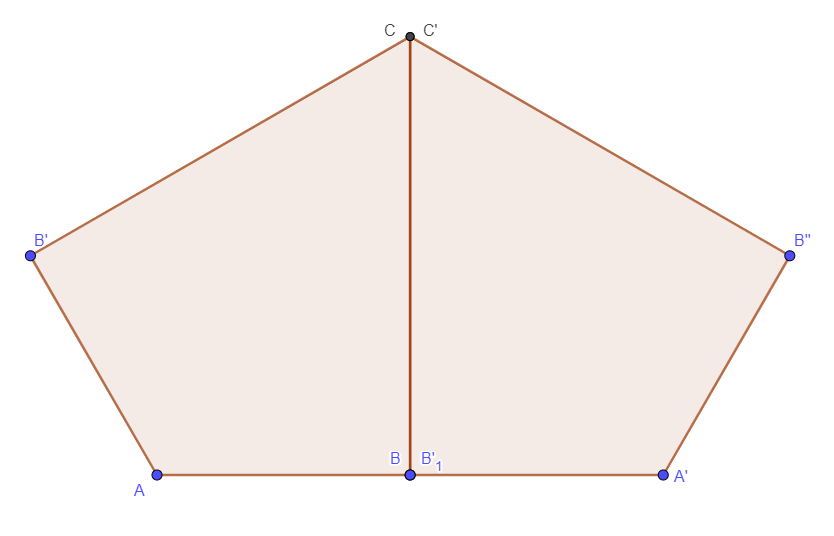}
    \caption{Primeira etapa na construção do Hat}
    \label{fig:hat_etapa1}
\end{figure}

\textbf{\emph{Hat}:} Basta fazermos reflexões sucessivas das pipas obtidas em relação aos seus lados. A primeira reflexão é da pipa construída anteriormente em relação ao lado $BC$ conforme ilustrado na Figura~\ref{fig:hat_etapa1}.

\begin{figure}[ht]
    \centering
    \includegraphics[width=0.55\textwidth]{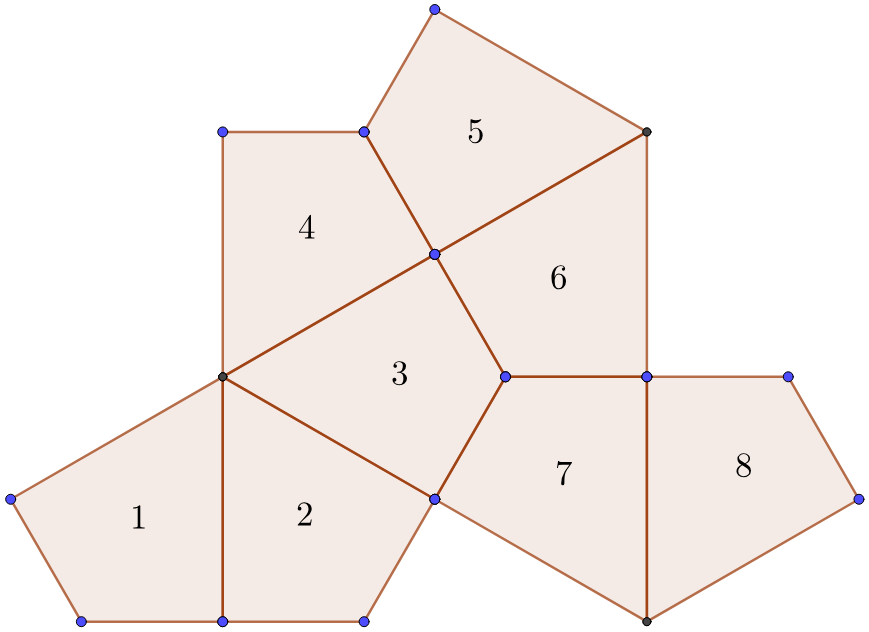}
    \caption{Segunda etapa na construção do \emph{Hat} (\href{https://www.geogebra.org/m/pthbn3du}{\emph{GeoGebra}})}
    \label{fig:hat_etapa2}
\end{figure}

Continuamos o processo até obtermos 8 pipas. Refletimos as pipas em relação aos lados conforme a ordem apresentada na Figura~\ref{fig:hat_etapa2}. Nós omitimos os rótulos de todos os elementos.

\begin{figure}[ht]
    \centering
    \includegraphics[width=0.5\textwidth]{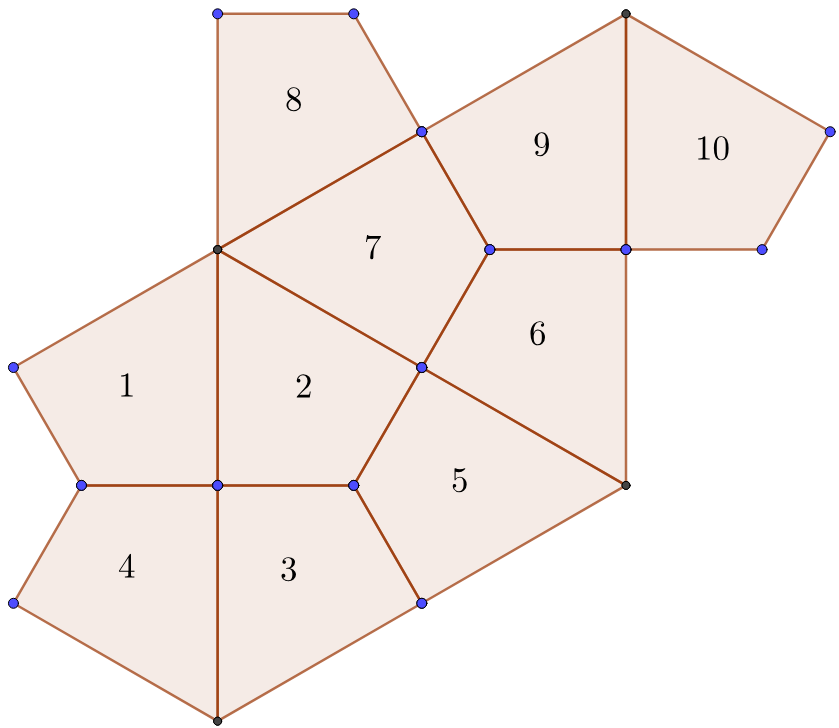}
    \caption{Construção do \emph{Turtle} (\href{https://www.geogebra.org/m/k5yrahfq}{\emph{GeoGebra}})}
    \label{fig:10pipa}
\end{figure}

\textbf{\emph{Turtle}:} Para construir o polígono $Tile(\sqrt{3},1)$, repetimos a construção da pipa de Laves (Figura~\ref{fig:pipa_etapa4}) e fazemos reflexões seguindo a ordem apresentada na Figura~\ref{fig:10pipa} até obtermos $10$ pipas.

\begin{figure}[ht]
    \centering
    \includegraphics[width=0.4\textwidth]{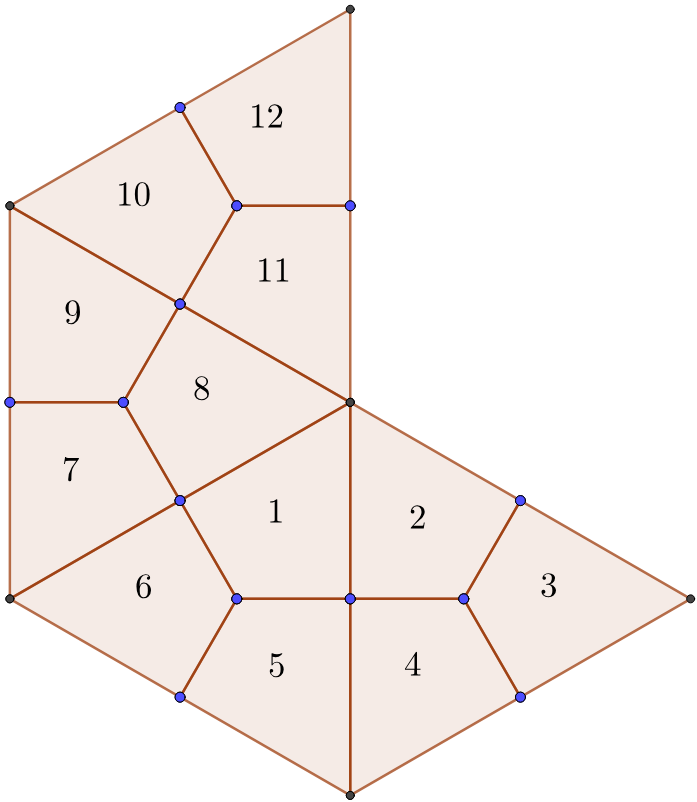}
    \caption{Ladrilho de Smith com $12$ pipas ou $Tile(0,1)$ (\href{https://www.geogebra.org/m/yfgzudu8}{\emph{GeoGebra}})}
    \label{fig:12pipa}
\end{figure}

\textbf{\emph{Ladrilho com $12$ pipas}:} Seguindo com as construções, vamos criar o ladrilho de Smith com 12 pipas ou $Tile(0,1)$. Novamente, começamos a construção com a pipa de Laves (Figura~\ref{fig:pipa_etapa4}) e seguimos fazendo reflexões. A ordem a ser seguida ao refletirmos as pipas é apresentada na Figura~\ref{fig:12pipa}.

\begin{figure}[ht]
    \centering
    \includegraphics[width=0.3\textwidth]{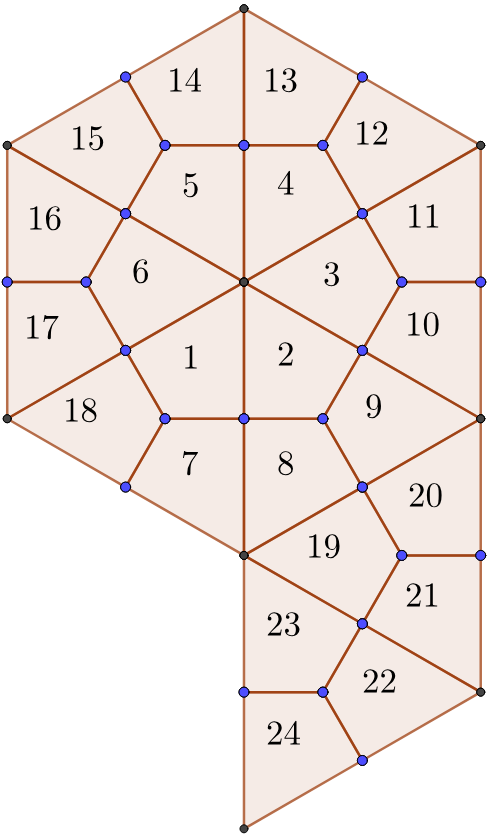}
    \caption{Ladrilho de Smith com $24$ pipas ou $Tile(1,0)$ (\href{https://www.geogebra.org/m/hmtmvdrk}{\emph{GeoGebra}})}
    \label{fig:24pipa}
\end{figure}

\textbf{Ladrilho com $24$ pipas:} Finalmente, iremos construir o ladrilho de Smith com 24 pipas ou $Tile(1,0)$. Como feito anteriormente, construímos a pipa de Laves (Figura~\ref{fig:pipa_etapa4}) e seguimos fazendo reflexões. A ordem a ser seguida ao refletirmos as pipas é apresentada na Figura~\ref{fig:24pipa}.

\section{Outros Polígonos de Smith}

Nesta seção, nós apresentaremos uma construção no GeoGebra que permite obter outros polígonos de Smith e seus colaboradores. Incrivelmente, estes pesquisadores descobriram uma infinidade de polígonos capazes de preencher o plano de forma aperiódica. Estes polígonos podem ser obtidos a partir da construção do polígono $Tile(a,1-a)$ com um parâmetro $a\in[0,1]$.

\begin{proposition}
O $Tile(a,b)$ é semelhante ao $Tile(ka,kb)$ para qualquer constante $k\neq 0$ \cite{smith2023aperiodic}.
\label{prop:tile_semelhantes}
\end{proposition}
\begin{proof}
De fato, os lados correspondentes dos dois polígonos $Tile(a,b)$ e $Tile(ka,kb)$ serão proporcionais com constante de proporcionalidade $k$.
\end{proof}

Queremos determinar uma constante $a\in [0,1]$ de modo que o polígono $Tile(a,1-a)$ corresponda ao \emph{Hat} ou ao polígono $Tile(1,\sqrt{3})$. Primeiramente, observamos que:
\begin{equation*}
    \text{O polígono } Tile(1,\sqrt{3}) \text{ é semelhante ao polígono } Tile(k,k\sqrt{3}) \text{ (Proposição~\ref{prop:tile_semelhantes})}
\end{equation*}
Vamos determinar uma constante $k$ tal que $a=k$ e $1-a=k\sqrt{3}$, desse modo:
\begin{equation}
    1-a=k\sqrt{3} \Rightarrow 1-k=k\sqrt{3} \Rightarrow 1=k+k\sqrt{3} \Rightarrow 1=k(1+\sqrt{3}) \Rightarrow k=\frac{1}{1+\sqrt{3}} \approx 0.37.
    \label{eq:k_hat}
\end{equation}
Em outras palavras, o \emph{Hat} (ou polígono $Tile(1,\sqrt{3})$) é semelhante ao polígono $Tile\left(\frac{1}{1+\sqrt{3}},\frac{\sqrt{3}}{1+\sqrt{3}}\right)$. Nós vamos usar esta propriedade do \emph{Hat} para apresentar uma outra forma de construção deste polígono. Antes disto, porém, observe que podemos mostrar de forma análoga que o \emph{Turtle} (ou polígono $Tile(\sqrt{3},1)$) é semelhante ao polígono $Tile\left(\frac{\sqrt{3}}{1+\sqrt{3}},\frac{1}{1+\sqrt{3}}\right)$. Mais ainda, o polígono $Tile(1,1)$ é semelhante ao polígono $Tile\left(\frac{1}{2},\frac{1}{2}\right)$. De fato, basta observarmos que o polígono $Tile(1,1)$ é semelhante ao polígono $Tile(k,k)$ e tentarmos escrever este polígono como $Tile(a,1-a)$. Ou seja, queremos uma constante $k$ tal que $a=k$ e $1-a=k$, então
\begin{equation}
    1-a=k \Rightarrow 1-k=k \Rightarrow 1=2k \Rightarrow k=\frac{1}{2}.
\end{equation}

\begin{figure}[ht]
    \centering
    \includegraphics[width=\textwidth]{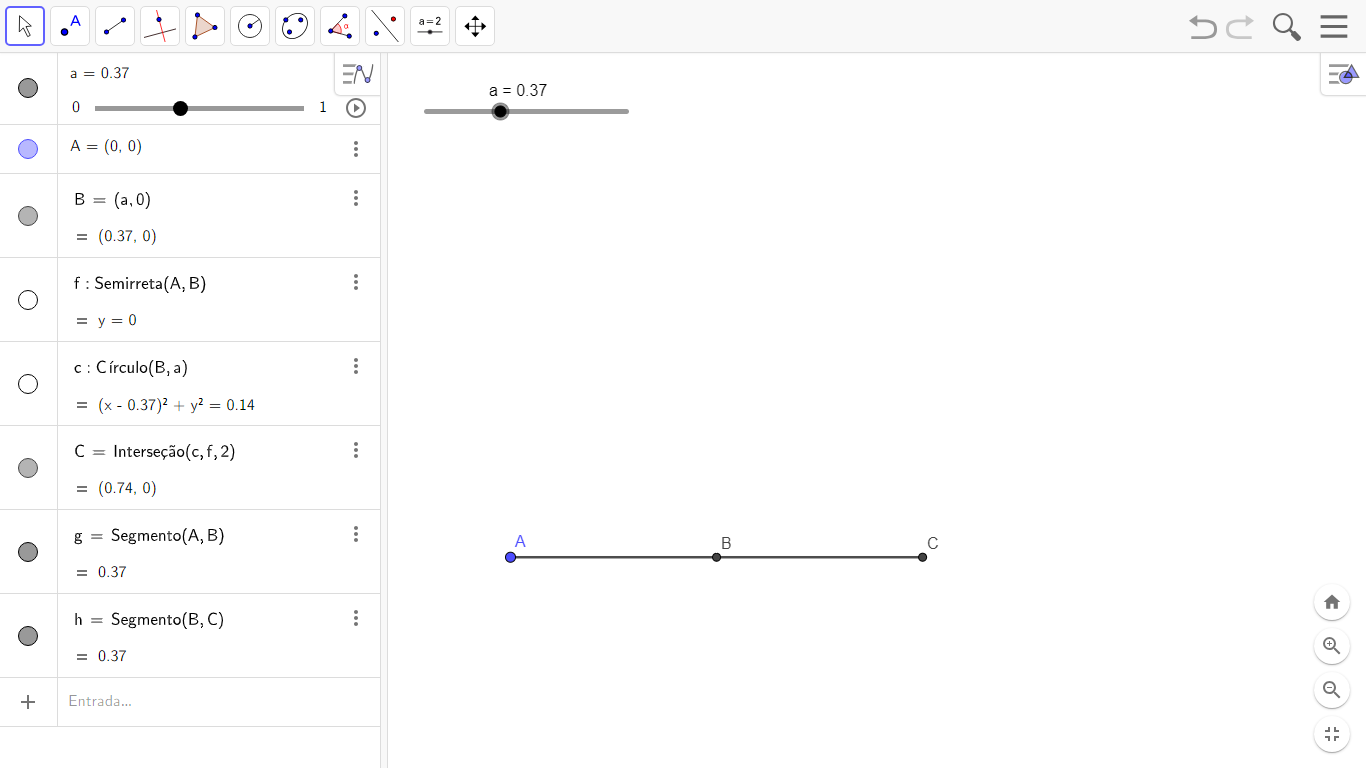}
    \caption{Primeira etapa da construção do $Tile(a,1-a)$}
    \label{fig:continuo_etapa1}
\end{figure}

\textbf{\emph{Hat e outros polígonos de Smith}:} Construímos inicialmente pontos $A$ e $B$ que distam $a$ um do outro e os unimos com um segmento. O valor de $a$ é definido por meio de um \emph{controle deslizante} com mínimo $0.001$ e máximo $0.999$. Usamos o valor inicial de $a=0.37$ como indicado na Equação~\eqref{eq:k_hat} a fim de obtermos o \emph{Hat} em um primeiro momento. Seguindo com a construção, traçamos a semirreta $\overrightarrow{AB}$, uma circunferência de centro $B$ e raio $a$, e determinamos o ponto $C$ de interseção destes dois objetos. Omitimos a semirreta e a circunferência e traçamos o segmento $BC$ (veja a Figura~\ref{fig:continuo_etapa1}).

\begin{figure}[ht]
    \centering
    \includegraphics[width=0.5\textwidth]{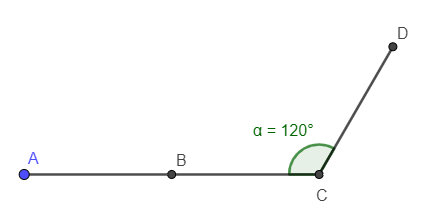}
    \caption{Segunda etapa da construção do $Tile(a,1-a)$}
    \label{fig:continuo_etapa2}
\end{figure}

Construímos então um ângulo de $120^{\circ}$ no sentido horário a partir dos pontos $B$ e $C$ (nesta ordem) obtendo assim um ponto $B'$. Em seguida, traçamos a semirreta $\overrightarrow{CB'}$ e uma circunferência de centro $C$ e raio $a$. Obtemos um ponto $D$ como a interseção da semirreta e da circunferência. Omitimos os objetos auxiliares criados e ligamos os pontos $C$ e $D$ com um segmento (veja a Figura~\ref{fig:continuo_etapa2}). Alternativamente, como $B'$ dista $a$ do ponto $C$, poderíamos simplesmente ter renomeado o ponto $B'$ chamando-o de $D$.

\begin{figure}[ht]
    \centering
    \includegraphics[width=0.65\textwidth]{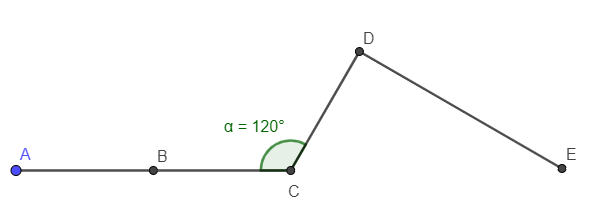}
    \caption{Terceira etapa da construção do $Tile(a,1-a)$}
    \label{fig:continuo_etapa3}
\end{figure}

Seguindo com a construção, traçamos uma reta perpendicular à $CD$ passando por $D$, uma circunferência de centro em $D$ e raio $1-a$. Obtemos o ponto $E$ de interseção destes objetos. Omitimos os objetos auxiliares e traçamos o segmento $DE$ (veja a Figura~\ref{fig:continuo_etapa3}).

\begin{figure}[ht]
    \centering
    \includegraphics[width=0.7\textwidth]{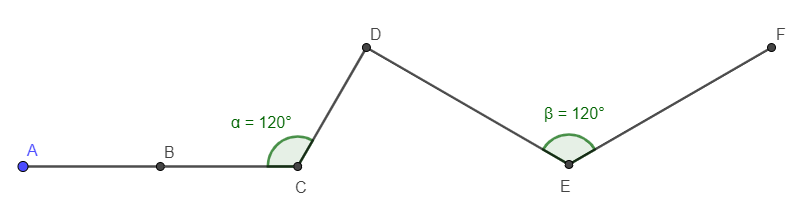}
    \caption{Quarta etapa da construção do $Tile(a,1-a)$}
    \label{fig:continuo_etapa4}
\end{figure}

Construímos um segundo ângulo de $120^{\circ}$ no senti horário a partir dos pontos $D$ e $E$ (nesta ordem) obtendo um ponto $D'$. Construímos então a semirreta $\overrightarrow{ED'}$, uma circunferência de centro $E$ e raio $1-a$ e determinamos a interseção destes dois, obtendo o ponto $F$. Finalizamos, esta etapa traçando o segmento $EF$ (veja a Figura~\ref{fig:continuo_etapa4}).

\begin{figure}[ht]
    \centering
    \includegraphics[width=0.7\textwidth]{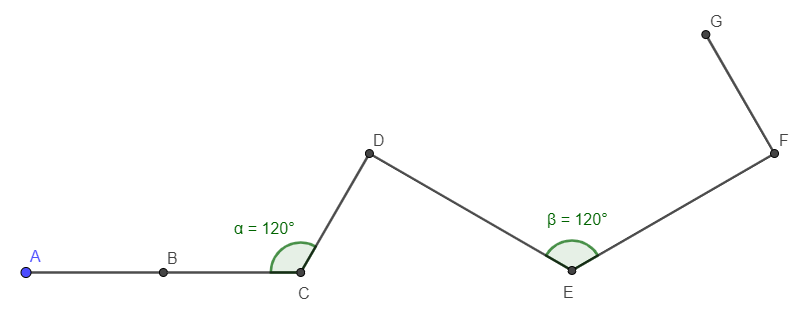}
    \caption{Quinta etapa da construção do $Tile(a,1-a)$}
    \label{fig:continuo_etapa5}
\end{figure}

Em seguida, traçamos uma reta perpendicular a $EF$ passando por $F$ e uma circunferência de centro $F$ e raio $a$ para determinar o ponto $G$ de interseção dos objetos. Construímos o segmento $FG$ (veja a Figura~\ref{fig:continuo_etapa5}).

\begin{figure}[ht]
    \centering
    \includegraphics[width=0.6\textwidth]{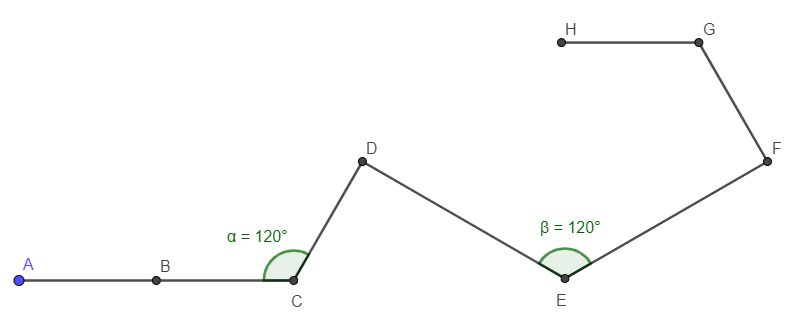}
    \caption{Sexta etapa da construção do $Tile(a,1-a)$}
    \label{fig:continuo_etapa6}
\end{figure}

Na sequência, vamos traçar uma paralela a $AB$ passando por $G$, uma circunferência de centro $G$ e raio $a$ e obter o ponto $H$ de interseção destes dois objetos. Construímos o segmento $GH$ (veja a Figura~\ref{fig:continuo_etapa6}).

\begin{figure}[ht]
    \centering
    \includegraphics[width=0.6\textwidth]{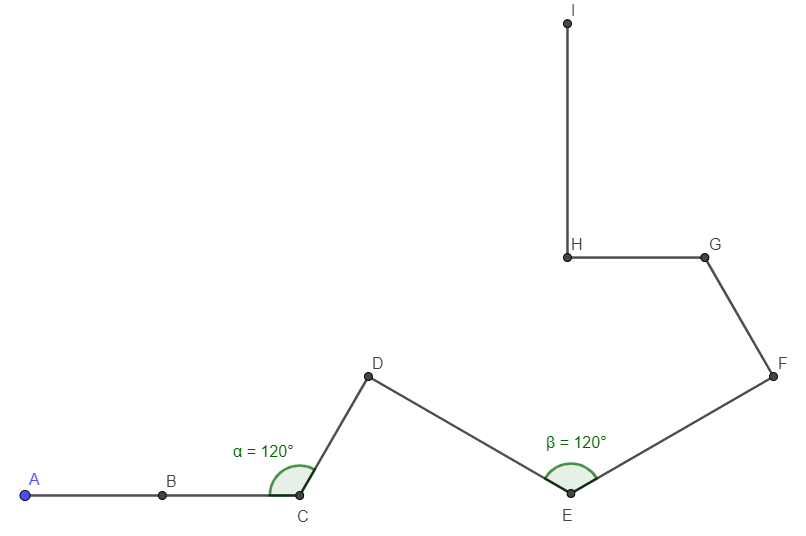}
    \caption{Sétima etapa da construção do $Tile(a,1-a)$}
    \label{fig:continuo_etapa7}
\end{figure}

Para continuar o processo, traçamos uma perpendicular a $GH$ por $H$ e uma circunferência de centro $H$ e raio $1-a$, obtemos assim o ponto $I$ de interseção destes objetos e construímos o segmento $HI$ (veja a Figura~\ref{fig:continuo_etapa7}).

\begin{figure}[h!]
    \centering
    \includegraphics[width=0.38\textwidth]{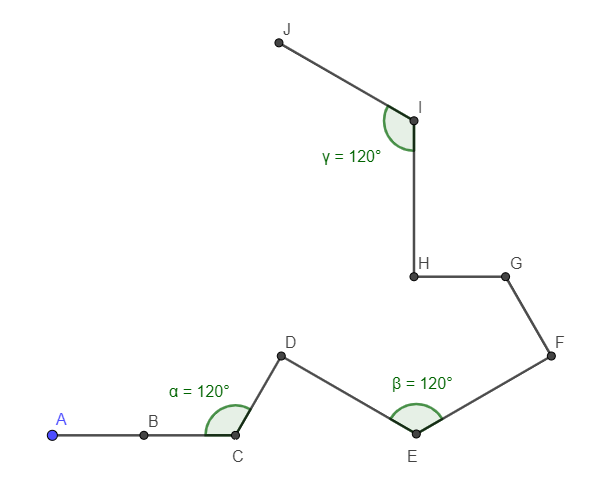}
    \caption{Oitava etapa da construção do $Tile(a,1-a)$}
    \label{fig:continuo_etapa8}
\end{figure}

Construímos então mais um ângulo de $120^{\circ}$ no senti horário a partir dos pontos $H$ e $I$ (nesta ordem) obtendo um ponto $J$. Traçamos o segmento $IJ$ (veja a Figura~\ref{fig:continuo_etapa8}).

\begin{figure}[ht]
    \centering
    \includegraphics[width=0.36\textwidth]{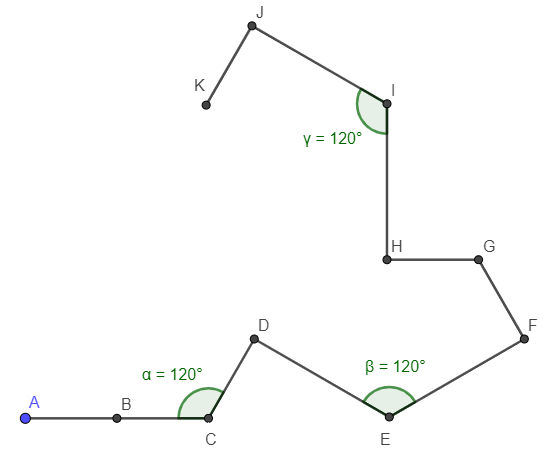}
    \caption{Nona etapa da construção do $Tile(a,1-a)$}
    \label{fig:continuo_etapa9}
\end{figure}

Traçamos uma reta perpendicular a $IJ$ por $J$, uma circunferência de cento $J$ e raio $a$ e obtemos o ponto $K$ de interseção destes objetos. Traçamos o segmento $JK$ (veja a Figura~\ref{fig:continuo_etapa9}).

\begin{figure}[ht]
    \centering
    \includegraphics[width=0.36\textwidth]{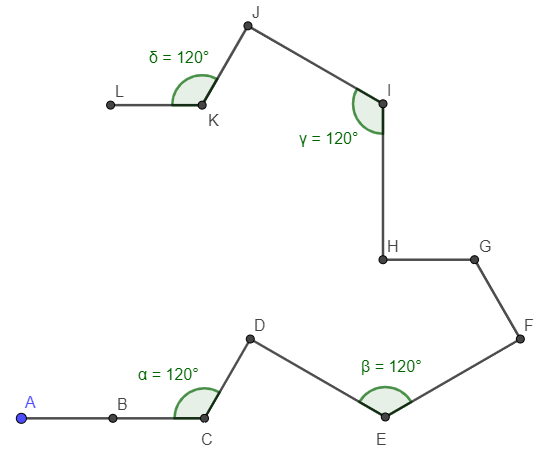}
    \caption{Décima etapa da construção do $Tile(a,1-a)$}
    \label{fig:continuo_etapa10}
\end{figure}

Construímos um ângulo de $120^{\circ}$ no senti anti-horário (ou $240^{\circ}$ no sentido horário) a partir dos pontos $J$ e $K$ (nesta ordem) obtendo um ponto $L$. Traçamos o segmento $KL$ (veja a Figura~\ref{fig:continuo_etapa10}).

\begin{figure}[h!]
    \centering
    \includegraphics[width=0.36\textwidth]{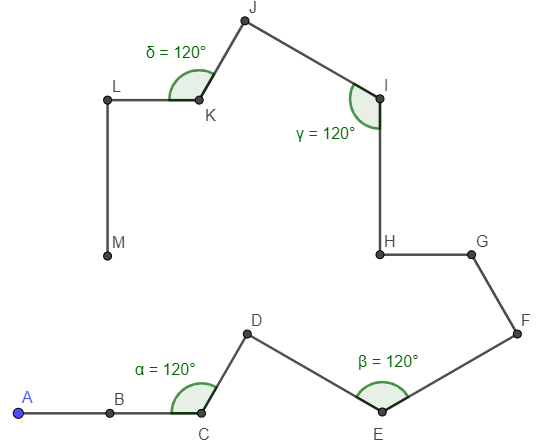}
    \caption{Décima primeira etapa da construção do $Tile(a,1-a)$}
    \label{fig:continuo_etapa11}
\end{figure}

Traçamos uma perpendicular a $KL$ por $L$, uma circunferência de centro $L$ e raio $1-a$ e definimos o ponto de interseção destes objetos como $M$. Traçamos o segmentos $LM$ (veja a Figura~\ref{fig:continuo_etapa11}).

\begin{figure}[ht]
    \centering
    \includegraphics[width=0.65\textwidth]{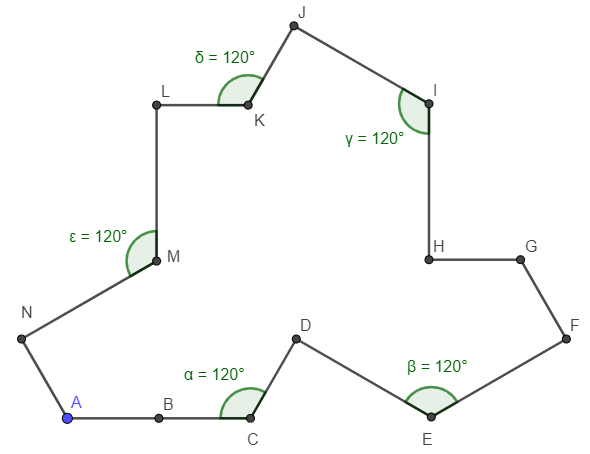}
    \caption{Décima segunda etapa da construção do $Tile(a,1-a)$ (\href{https://www.geogebra.org/m/tuqybfrb}{\emph{GeoGebra}})}
    \label{fig:continuo_etapa12}
\end{figure}

Construímos um ângulo de $120^{\circ}$ no senti anti-horário (ou $240^{\circ}$ no sentido horário) a partir dos pontos $L$ e $M$ (nesta ordem) obtendo um ponto $N$. Traçamos o segmento $MN$ (veja a Figura~\ref{fig:continuo_etapa12}). E, finalmente, traçamos o segmento $NA$.

Com o passo a passo descrito nesta construção obtemos o polígono $Tile(1,\sqrt{3})$ ou ladrilho com $8$ pipas de Laves. Mas, modificando dinamicamente o parâmetro $a$ obtemos um conjunto de polígonos que são capazes de preencher o plano de forma unicamente aperiódica\footnote{Vídeo resultante da variação do parâmetro $a$: \url{https://youtu.be/99d9R7PxgwI}}. Em outras palavras, os polígonos de Smith e seus colaboradores podem ser construídos a partir de uma outra perspectiva, diferente daquela fundamentada nas pipas de Laves.

\begin{figure}[ht]
    \centering
    \begin{minipage}[b]{0.3\textwidth}
    \centering
    \includegraphics[width=0.77\textwidth]{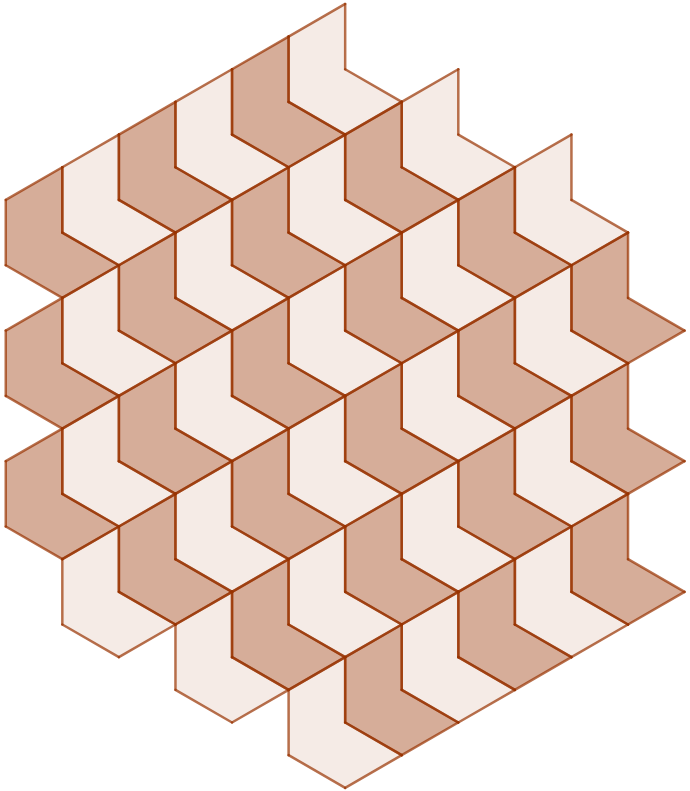}
    \caption{$Ladrilho(0,1)$ em uma tesselação periódica}
    \label{fig:tesselacao_12pipa_geoebra}
    \end{minipage} \quad
    \begin{minipage}[b]{0.3\textwidth}
    \centering
    \includegraphics[width=0.9\textwidth]{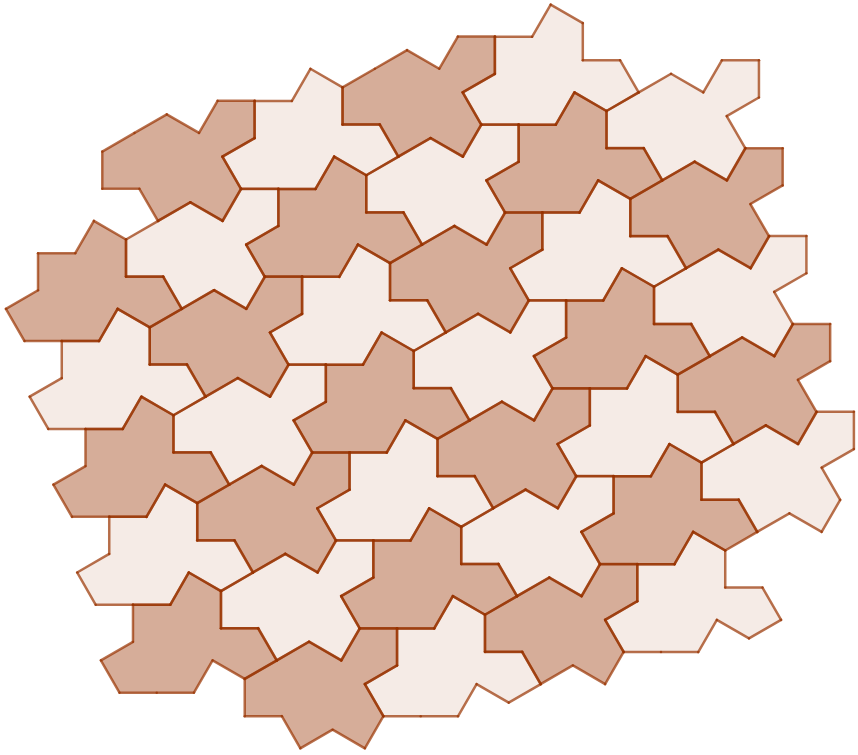}
    \caption{$Ladrilho(1,1)$ em uma tesselação periódica}
    \label{fig:tesselacao_10pipa_geoebra}
    \end{minipage} \quad
    \begin{minipage}[b]{0.3\textwidth}
    \centering
    \includegraphics[width=0.8\textwidth]{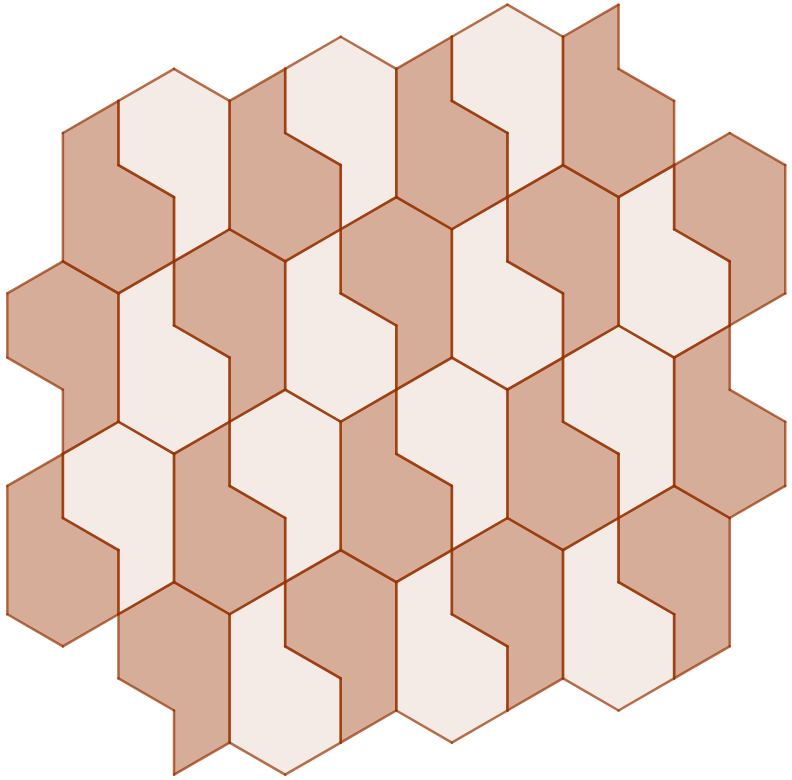}
    \caption{$Ladrilho(1,0)$ em uma tesselação periódica}
    \label{fig:tesselacao_24pipa_geoebra}
    \end{minipage}
\end{figure}

Cada polígono obtido a partir de uma constante $a\in [0,1]$ é capaz de gerar uma tesselação unicamente aperiódica do plano com exceção dos casos $Tile(1,0)$, $Tile(0,1)$ e $Tile(1,1)$ que também admitem tesselações periódicas como indicado nas Figuras~\ref{fig:tesselacao_12pipa_geoebra}, \ref{fig:tesselacao_10pipa_geoebra} e \ref{fig:tesselacao_24pipa_geoebra}.

\section{Conclusão}

Neste artigo, nós apresentamos em detalhes algumas formas de construção dos polígonos de Smith e seus colaboradores. A primeira construção feita através de uma pipa em particular é bastante natural quando consideramos a tesselação dual de Laves da tesselação semirregular $(3.4.6.4)$. Os famosos polígonos Hat e Turtle podem ser construídos respectivamente com $8$ e $10$ pipas seguindo ordens específicas de construção. A segunda construção fundamentada nos lados de cada um dos polígonos permite com a mudança de um parâmetro obter diversos polígonos que admitem apenas uma tesselação aperiódica do plano. Neste caso, nós mostramos quais são os valores de parâmetro que dão origem ao Hat e ao Turtle. Ao longo do texto as construções foram feitas passo a passo no software GeoGebra o que permite a reprodução das construções por outros professores e pesquisadores. Além disso, a dinâmica fornecida pelo programa ajuda a visualizar a classe de polígonos descobertos por Smith e seus colaboradores com a mudança de um simples parâmetro.

\section*{Agradecimentos}
Os autores desse artigo gostariam de agradecer ao Programa de Educação Tutorial - FNDE.

\end{document}